\documentclass[a4paper]{amsart}

\author{Samuel Boissi{\`e}re}

\title{Towards the multiplicative behavior of the K-theoretical McKay correspondence}

\address{Samuel Boissi\`{e}re, Fachbereich f\"{u}r Mathematik, Staudinger Weg $9$,
Johannes Gutenberg-Universit\"{a}t Mainz, $55099$ Mainz, Germany}

\email{boissiere@mathematik.uni-mainz.de}

\urladdr{http://sokrates.mathematik.uni.mainz.de/$\sim$samuel}

\usepackage{amsmath}
\usepackage[all]{xy}
\usepackage[english]{babel}
\usepackage{mathrsfs}
\usepackage{eufrak}
\usepackage{youngtab}
\usepackage{rotating}

\DeclareMathOperator{\Hilb}{Hilb} \DeclareMathOperator{\Coeff}{Coeff}
\DeclareMathOperator{\codim}{codim} \DeclareMathOperator{\Supp}{Supp}
\DeclareMathOperator{\Hom}{Hom} \DeclareMathOperator{\Spec}{Spec}
 \DeclareMathOperator{\Res}{Res}
\DeclareMathOperator{\End}{End}
\DeclareMathOperator{\age}{age}\DeclareMathOperator{\GL}{GL}
\DeclareMathOperator{\Id}{Id} \DeclareMathOperator{\gr}{gr}
\DeclareMathOperator{\Sp}{Sp}

\newcommand{\IC}{\mathbb{C}}
\newcommand{\IN}{\mathbb{N}}
\newcommand{\IQ}{\mathbb{Q}}
\newcommand{\IZ}{\mathbb{Z}}

\newcommand{\kg}{\mathfrak{g}}
\newcommand{\kh}{\mathfrak{h}}
\newcommand{\kU}{\mathfrak{U}}

\newcommand{\hilb}{\Hilb^n(\IC^2)}
\newcommand{\Snhilb}{S_n\text{-}\Hilb\IC^{2n}}
\newcommand{\LambdaQ}{\Lambda_{\IQ(q,t)}}
\newcommand{\eins}{\mathbf{1}}

\newcommand{\tH}{\widetilde{H}}
\newcommand{\tK}{\widetilde{K}}

\newcommand{\cC}{\mathcal{C}}
\newcommand{\cF}{\mathcal{F}}
\newcommand{\cO}{\mathcal{O}}

\newcommand{\opD}{\mathscr{D}}
\newcommand{\opE}{\mathscr{E}}

\theoremstyle{plain}
\newtheorem{theorem}{Theorem}[section]

\newtheorem{proposition}[theorem]{Proposition}
\newtheorem{lemma}[theorem]{Lemma}
\newtheorem{remark}[theorem]{Remark}

\newtheorem{conjecture}[theorem]{Conjecture}

\begin{document}

\begin{abstract}
When the quotient of a symplectic vector space by the action of a finite
subgroup of symplectic automorphisms admits as a crepant projective resolution
of singularities the Hilbert scheme of regular orbits of Nakamura, then there
is a natural isomorphism between the Grothendieck group of this resolution and
the representation ring of the group, given by the Bridgeland-King-Reid map.
However, this isomorphism is not compatible with the ring structures. For the
Hilbert scheme of points on the affine plane, we study the multiplicative
behavior of this map.
\end{abstract}



\maketitle

\pagestyle{myheadings}

\markboth{SAMUEL BOISSI\`ERE}{Multiplicative behavior of the K-theoretical
McKay correspondence}

\section{Introduction}

Let $V$ be a finite-dimensional complex symplectic vector space and $G\subset
\Sp(V)$ a finite subgroup of automorphisms respecting the symplectic form. The
quotient $V/G=\Spec \cO(V)^G$ is usually a singular variety. The
Bridgeland-King-Reid theorem (\cite{BKR}) asserts that in some special cases,
one can take as a crepant and projective resolution of singularities the
\emph{Hilbert scheme of $G$-regular orbits} $Y:=G\text{-}\Hilb V$ defined by
Nakamura (\cite{Nm}) and there exists a natural Fourier-Mukai functor inducing
an isomorphism between the Grothendieck group $K(Y)$ of coherent sheaves on $Y$
and the group of representations $R(G)$. This induces a $\IQ$-linear
isomorphism $K(Y)\underset{\IZ}{\otimes}\IQ \cong
R(G)\underset{\IZ}{\otimes}\IQ$. In the sequel, we always work with \emph{rational} coefficients.

Nakajima (\cite[Question $4.23$]{N1}) formulates the following question:

\noindent\textbf{Ring structure problem.}  The Grothendieck group $K(Y)$ has a
natural ring structure induced by the tensor product of vector bundles on $Y$,
whereas the ring structure on $R(G)$ is given by the tensor product of
representations. The preceding isomorphism is not compatible with these ring
structures (see also \cite[Problem $1.5$]{GK}). How to describe the product
induced on $R(G)$ by this map?

We study this question in the particular case $V=\IC^n\otimes \IC^2$ with the
permutation action of the symmetric group $G=S_n$ on the first factor and the canonical symplectic
structure coming from the second factor. The Hilbert scheme $Y:=\hilb$ provides a natural symplectic
resolution of singularities, isomorphic to the Hilbert scheme of regular orbits
$\Snhilb$, providing concrete informations about the Bridgeland-King-Reid map
(see Haiman \cite{H6,H5}). In this context, we replace the spaces $R(S_n)$ by the \emph{space of symmetric functions}
$\Lambda\cong\IQ[p_1,p_2,\ldots]$ and we work in the extended ring $\LambdaQ$ equipped
with the basis of \emph{modified Macdonald polynomials} $\tH_\mu$. In order to
study the problem, we consider a bundle $F$ on $\hilb$ and the induced operator
on $K(\hilb)$ given by tensor product $G\mapsto F\otimes G$. We look at a
description of the operator $\opE_F$ one gets on $\Lambda^n$ through the
isomorphism $\Theta:\Lambda^n\rightarrow K(\hilb)$ induced by the
Bridgeland-King-Reid map:
$$
\xymatrix{\Lambda^n \ar[rr]^-\Theta\ar[d]_{\opE_F}& & K(\hilb)\ar[d]^{F\otimes - } \\
\Lambda^n \ar[rr]^-\Theta && K(\hilb)}
$$
There is a natural action of the torus $(\IC^*)^2$ on $\hilb$ with a finite
number of fixed points indexed by \emph{partitions} $\mu$ of $n$. The vector
space $K(\hilb)$ has a basis of $(\IC^*)^2$-linearized bundles (\S
\ref{ss:Hilbert orbits}), so it is enough to work with a $(\IC^*)^2$-linearized
vector bundle $F$ on $\hilb$. Denote by $F_\mu\in \IQ[q^{\pm 1},t^{\pm 1}]$ the
representation associated to the action on the fibre at each fixed point and
define an operator $\nabla_F:\LambdaQ^n\rightarrow \LambdaQ^n$ by:
$\nabla_F \tH_\mu=F_\mu \tH_\mu$, $\forall \mu\vdash n$.
Denote by $\nabla^*_F$ the adjoint operator of $\nabla_F$ for the scalar product
$\langle,\rangle$ on $\Lambda$ and by $\omega$ the natural involution on
$\Lambda$ (see \S\ref{ss:Ring Symmetric}). We get the following description:

\begin{theorem} The operator $\opE_F$ defined by:
$$
F\otimes \Theta(y)=\Theta(\opE_F(y))\quad \forall y\in \Lambda^n
$$
in given by the following formula:
$$
\opE_F:=\left(\omega \nabla_F^*\omega\right)_{q=1,t=1}.
$$
\end{theorem}

In particular, for the tautological bundle $B_n$ on $\hilb$ we find a
differential operator for all $n$ together, which appears as a K-theoretic
analogue of a result of Lehn on the cohomological action of the Chern character
of $B_n$ (see theorem \ref{th:operatorD} and the operator $\opD$):
\begin{theorem}
Consider the following differential operator on $\Lambda$:
$$
\opE:=\Coeff\left(t^0,\left(\sum\limits_{r\geq 1}
rp_rt^r\right)\exp\left(\sum\limits_{r\geq 1}\frac{\partial}{\partial p_r}
t^{-r}\right)\right).
$$
Then for any $y\in \Lambda^n$ one has:
$$
B_n\otimes\Theta(y)=\Theta(\opE(y)).
$$
\end{theorem}

In the last sections of this paper, we investigate  the relations between our
operator $\opE$ and the operator $\opD$ of Lehn in relation with the McKay
correspondences and the Chern character, as pictured in the following diagram:
$$
\xymatrix{ \Lambda^n \ar@(l,u)[]|{\opD} \ar[rr]^-\Psi && H^*(\hilb) \\
\Lambda^n \ar@(d,l)[]|{\opE} \ar@{.>}[u]^\Gamma \ar[rr]^-\Theta &&  K(\hilb)
\ar[u]^{ch} }
$$
By making use of our preceding results \cite{SB2} about the induced map
$\Gamma$, we give  some other precise informations about the ring structure
induced on $\Lambda^n$ by $\Theta$, as such as numerical tables in the basis of
Newton and Schur functions.

I thank C. Sorger and M. Lehn for their continuous encouragement and many helpful
discussions, and the anonymous referee for his carefully reading and his suggestions.

\section{Symmetric functions}
\label{s:symmetric}

\subsection{The ring of symmetric functions} (\cite{McDo1,M})
\label{ss:Ring Symmetric}

Take independent indeterminates $x_1,\ldots,x_r$, let the symmetric group $S_r$
act by permutation on $\IQ[x_1,\ldots,x_r]$ and denote the invariant ring by
$\Lambda_r:=~\IQ[x_1,\ldots,x_r]^{S_r}$. This ring is naturally graded by
degree and we denote by $\Lambda^n_r$ the vector subspace generated by degree
$n$ homogeneous symmetric polynomials. By adjoining other indeterminates one
constructs the projective limit $\Lambda^n:=\underleftarrow{\lim}\Lambda^n_r$.
Then the \emph{ring of symmetric functions} is defined as
$\Lambda:=\bigoplus\limits_{n\geq 0}\Lambda^n$.

A \emph{partition} of an integer $n$ is a non-increasing sequence of non-negative
integers $\lambda:=(\lambda_1,\ldots,\lambda_k)$ such that
$\sum\limits_{i=1}^k\lambda_i=n$ (we write $\lambda\vdash n$). The $\lambda_i$
are the \emph{parts} of the partition. If necessary, we extend a partition with
zero parts. The number $l(\lambda)$ of non-zero parts is the \emph{length} of
the partition and the sum $|\lambda|$ of the parts is the \emph{weight}. If a
partition $\lambda$ has $\alpha_1$ parts equal to $1$, $\alpha_2$ parts equal
to $2$, $\ldots$ we shall also denote it by
$\lambda:=(1^{\alpha_1},2^{\alpha_2},\ldots)$.

The \emph{Young diagram} of a partition $\lambda$ is defined by:
$$
D(\lambda):=\{(i,j) \in \IN\times\IN \,|\, j<\lambda_{i+1} \}.
$$
In the representation of such a diagram, we follow a matrix convention:
$$
\Yvcentermath1\young(~xaa,~l~,~)
\begin{array}{ccc}
\lambda=(4,3,1) && x=(0,1)\\
|\lambda|=8 && a(x)=2\quad l(x)=1 \\
l(\lambda)=3 && h(x)=4
\end{array}
$$
For each cell $x\in D(\lambda)$, the \emph{arm} $a(x)$ is the number of cells
on the right of $x$, the \emph{leg} $l(x)$ is the number of cells below $x$,
the \emph{hook length} is $h(x):=1+a(x)+l(x)$. We shall also make use of the
number $n(\lambda):=\sum\limits_{i\geq 1}(i-1)\lambda_i$.

Set $p_0=1$ and define for $k\geq 1$ the \emph{power sum}
$p_k:=\sum\limits_{i\geq 1} x_i^k$. For a partition
$\lambda=(\lambda_1,\ldots,\lambda_k)$, the \emph{Newton function} is the
product $p_\lambda:=p_{\lambda_1}\cdots p_{\lambda_k}\in \Lambda^{|\lambda|}$.
The Newton functions form a basis of $\Lambda$ and
$\Lambda\cong\IQ[p_1,p_2,\ldots]$. We denote by $\omega:\Lambda\rightarrow\Lambda$
the algebra isomorphism characterized by $\omega p_k=-p_k$. Another natural
basis of $\Lambda$ indexed by partitions is given by the \emph{Schur functions}
$s_\lambda$.

For a partition $\lambda=(1^{\alpha_1},2^{\alpha_2},\ldots)$, set
$z_\lambda:=\prod\limits_{r\geq 1} \alpha_r!r^{\alpha_r}$ and define a scalar
product on $\Lambda$ by $\langle p_\lambda,p_\mu\rangle=\delta_{\lambda,\mu}
z_\lambda$ where $\delta_{\lambda,\mu}$ is the Kronecker symbol. Then the basis
of Schur functions is orthonormal.

Let $\cC(S_n)$ be the $\IQ$-vector space of \emph{class functions} on $S_n$.
Since conjugacy classes in $S_n$ are indexed by partitions, the functions
$\chi_\lambda$ taking the value $1$ on the conjugacy class $\lambda$ and otherwise $0$
form a basis of $\cC(S_n)$. Let $R(S_n)$ be the $\IQ$-algebra of
representations of $S_n$, with the product induced by the tensor product of
representations. By associating to each representation of $S_n$ its
\emph{character} we get an isomorphism $\chi:R(S_n)\rightarrow \cC(S_n)$. The
\emph{Frobenius morphism} is the isomorphism $\Phi:\cC(S_n)\rightarrow
\Lambda^n$ characterized by $\Phi(\chi_\lambda)=z_\lambda^{-1}p_\lambda$. Let
$\chi^\lambda$ be the class function such that $\Phi(\chi^\lambda)=s_\lambda$
and $\chi^\lambda_\mu$ the value of $\chi^\lambda$ at the conjugacy class
$\mu$. The representations $V^\lambda$ of character $\chi^\lambda$ are the
irreducible representations of $S_n$ and we have the following base-change
formulas:
\begin{align*}
p_\mu&=\sum_{\lambda \vdash n} \chi^\lambda_\mu s_\lambda \text{ (Frobenius formula}), \\
s_\lambda&=\sum_{\mu \vdash n} z_\mu^{-1} \chi^\lambda_\mu p_\mu \text{
(inverse Frobenius formula)}.
\end{align*}

\subsection{Plethystic substitutions} (\cite{H1,McDo1})
\label{ss:Plethystic}

The identification $\Lambda=\IQ[p_1,p_2,\ldots]$ allows to specialize the
$p_k$'s to elements of any $\IQ$-algebra: the specialization extends uniquely
to an algebra homomorphism on $\Lambda$. For a formal Laurent series $E$ in
indeterminates $t_1,t_2,\ldots$ we define $p_k[E]$ to be the result of
replacing each indeterminate $t_i$ by $t_i^k$. Extending the specialization to
any symmetric function $f\in \Lambda$, we obtain the \emph{plethystic
substitution} of $E$ in $f$, denoted by $f[E]$. Our convention is that in a
plethystic substitution, $X$ stands for the sum of the original indeterminates
$x_1+x_2+\cdots$, so that $p_k[X]=p_k$.

\subsection{Macdonald polynomials}(\cite{H1,McDo2})
\label{ss:Macdonald}

We introduce indeterminates $q,t$ and consider the ring
$\Lambda_{\IQ(q,t)}:=\Lambda\underset{\IQ}{\otimes} \IQ(q,t)$. The scalar
product and the plethystic substitutions naturally extend to this setup. For a
partition $\mu$, we define $B_\mu(q,t):=\sum\limits_{(i,j)\in D(\mu)}t^iq^j$.
Set:
$$
\Omega:=\exp\left(\sum\limits_{k\geq 1} \frac{p_k}{k}\right)
$$
and define a linear operator $\Delta:\Lambda_{\IQ(q,t)}\rightarrow
\Lambda_{\IQ(q,t)}$ by:
$$
\Delta f=\left. f \left[ X+\frac{(1-t)(1-q)}{z} \right] \Omega[-zX]
\right|_{z^0}.
$$
The \emph{modified Macdonald polynomial} $\tH_\mu$ is the eigenvector of
$\Delta$ corresponding to the eigenvalue $1-(1-q)(1-t)B_\mu(q,t)$. These
polynomials form a basis of $\Lambda_{\IQ(q,t)}$ and decompose in the basis of
Schur functions as:
$$
\tH_\mu=\sum_{\lambda \vdash n} \tK_{\lambda,\mu}s_\lambda,
$$
where the $\tK_{\lambda,\mu}\in \IN[q,t]$ are called $q,t$-\emph{Kostka
polynomials}.

The modified Macdonald polynomials satisfy the following duality formula
(\cite[\S$3.5$ formula ($53$)]{H2}):
\begin{equation}\label{eq:duality}
\langle \tH_\mu,\omega
\tH_\lambda[X(1-q)(1-t)]\rangle=\delta_{\lambda,\mu} a_\lambda
\end{equation}
for some non-zero coefficients $a_\lambda\in \IQ(q,t)$.

\section{Hilbert schemes on the affine plane}
\label{s:Hilbert}

\subsection{Hilbert scheme of points}(\cite{B,EG,ES,Fo,G,LS,N1,N3})
\label{ss:Hilbert points}

The \emph{Hilbert scheme of $n$ points in the affine plane} $\hilb$ is the
smooth quasi-projective manifold of complex dimension $2n$ parameterizing
length $n$ finite subschemes on the plane $\IC^2$. The first projection
$B_n:=pr_{1*}\cO_{\Xi_n}$ of the universal family $\Xi_n\subset \hilb\times
\IC^2$ is the rank $n$ \emph{usual tautological bundle} on $\hilb$.

The manifold $\hilb$ has no odd singular cohomology; its even cohomology has no
torsion and is generated by algebraic cycles. We denote by $H^*(\hilb)$ the
cohomology ring and by $K(\hilb)$ the Grothendieck algebra of algebraic vector
bundles (or equivalently of coherent sheaves), both with \emph{rational} coefficients.
Denote by $ch:K(\hilb)\overset{\sim}{\longrightarrow} H^*(\hilb)$ the Chern
character.

Denote by $S^n\IC^2:=\IC^{2n}/S_n$ the Mumford quotient parameterizing length $n$
zero-cycles in $\IC^2$. The \emph{Hilbert-Chow} morphism $\rho:\hilb\rightarrow
S^n\IC^2$ is a symplectic resolution of singularities.

There is a natural isomorphism
$$
\Psi:\Lambda^n \longrightarrow H^*(\hilb)
$$
constructed by use of geometric operators acting on the total sum of cohomology
of Hilbert schemes (see Nakajima \cite{N3}). In order to describe the ring
structure of $H^*(\hilb)$, Lehn-Sorger \cite{LS} use the following operator:

\begin{theorem}[Lehn]\emph{(\cite[Theorem $4.10$]{L})}
\label{th:operatorD} Consider the following differential operator on $\Lambda$:
$$
\opD:=\Coeff\left(t^0,-\left(\sum\limits_{r\geq 1}
p_rt^r\right)\exp\left(-\sum\limits_{r\geq 1}r\frac{\partial}{\partial p_r}
t^{-r}\right)\right).
$$
Then for any $y\in \Lambda^n$ one has:
$$
ch(B_n)\cup\Psi(y)=\Psi(\opD(y)).
$$
\end{theorem}

As explained in \cite[Theorem $4.1$ \& \S 6]{LS}, this gives a complete
description of the ring structure with generators and relations, since the
homogeneous components $ch_k(B_n)$ of the Chern character of $B_n$ generate the
ring $H^*(\hilb)$.

\subsection{Hilbert scheme of regular orbits} (\cite{BKR,H6,INj,IN,Nm})
\label{ss:Hilbert orbits}

The \emph{Hilbert scheme of $S_n$-regular orbits} $\Snhilb$ is defined as the
closure in the Hilbert scheme $\Hilb^{n!}(\IC^{2n})$ of the open set of
$S_n$-free orbits and is isomorphic to the Hilbert scheme $\hilb$ (see
\cite[Theorem $5.1$]{H6}). Denote the universal family by $Z_n\subset
\Snhilb\times \IC^{2n}$ and set $P_n:=p_*\cO_{Z_n}$ considered as the rank $n!$
\emph{unusual tautological bundle} on $\hilb$. This bundle is equipped with a
natural $S_n$-action inducing the regular representation on each fiber.
Consider the diagram:
$$
\xymatrix{Z_n\ar[r]^-q \ar[d]_-p&\IC^{2n} \ar[d]\\
\hilb\ar[r]^-\rho&S^n\IC^2}
$$
Denote by $D^b(\hilb)$ the derived category of coherent sheaves and by
$D^b_{S_n}(\IC^{2n})$ the derived category of $S_n$-equivariant coherent
sheaves. In this situation, we can apply the Bridgeland-King-Reid theorem
(\cite{BKR}) and get an isomorphism of Grothendieck groups
$$
\Upsilon:=q_!\circ p^!:K(\hilb)\rightarrow K_{S_n}(\IC^{2n}).
$$

Identifying $K_{S_n}(\IC^{2n})$ with the Grothendieck group of $S_n$-equivariant 
$\cO(\IC^{2n})$-modules of finite type one has (see
\cite{H5}):
$$
\Upsilon(F)=\sum\limits_{i\geq 0} (-1)^i H^i(\hilb,P_n\otimes F).
$$

Consider the following composition of vector space isomorphisms:
$$
\Theta:\Lambda^n\xrightarrow{\Phi^{-1}} \cC(S_n) \xrightarrow{\chi^{-1}} R(S_n)
\xrightarrow{\tau^{-1}} K_{S_n}(\IC^{2n})\xrightarrow{\Upsilon^{-1}} K(\hilb),
$$
where $\tau$ is the Thom isomorphism (here it is the restriction to a fibre,
see \cite[Theorem $5.4.17$]{CG}). The $S_n$-action on $P_n$ induces an
isotypical decomposition:
$$
P_n=\bigoplus_{\mu\vdash n} \mathbf{V}^\mu\otimes P_\mu,
$$
where $\mathbf{V}^\mu$ is the trivial bundle with fibre $V^\mu$ on $\hilb$ and
$P_\mu:=\Hom_{S_n}(\mathbf{V}^\mu,P_n)$. Then as in Ito-Nakajima \cite[Formula
($5.3$)]{IN} one has:

\begin{proposition}\label{prop:IN} For $\mu\vdash n$, $\Theta(s_\mu)=P_\mu^*$.
In particular, the dual bundles $P_\mu^*$ form a basis of $K(\hilb)$.
\end{proposition}

\subsection{Torus action on the Hilbert scheme of points}(\cite{ES,H5})
\label{ss:Torus action}

The torus $T^2:=(\IC^*)^2$ acts on $\IC[x,y]$ by $(t,q).x=tx, (t,q).y=qy$ for
$(t,q)\in~T^2$. It induces a natural action on $\hilb$ with finitely many fixed
points $\xi_\lambda$ parameterized by the partitions $\lambda$ of $n$. The
action extends to all natural objects at issue over $\IC^2$. In particular, the
isomorphism $\Upsilon$ is compatible with this action and induces an
isomorphism $\Upsilon:K_{T^2}(\hilb)\rightarrow K_{S_n\times T^2}(A_n)$ where
$A_n:=\cO(\IC^{2n})$.

Let $F$ be a $T^2$-linearized vector bundle on $\hilb$. Each fibre
$F(\xi_\lambda)$ has a structure of representation of $T^2$ and we denote by
$F_\lambda$ the polynomial in $q,t$ corres\-ponding to $F(\xi_\lambda)$ in the
identification $R(T^2)\cong \IQ[q^{\pm 1},t^{\pm 1}]$.

\subsection{Character formulas}(\cite{H5})

There is a generalization of the Frobenius isomorphism for $S_n\times
T^2$-equivariant $A_n$-modules of finite type. Let $M$ be a $S_n\times
T^2$-equivariant $A_n$-module of finite type and
$M=\bigoplus\limits_{i,j}M_{i,j}$ its isotypical decomposition with respect to
$T^2$: each $M_{i,j}$ is a finite-dimensional representation of $S_n$ where
$T^2$ acts by multiplication by $t^iq^j$. The \emph{Frobenius series} of $M$ is
the formal series:
$$
\cF_M:=\sum\limits_{i,j} t^iq^j \Phi(M_{i,j})\in \Lambda\otimes \IQ[[q^{\pm
1},t^{\pm 1}]].
$$
We shall use the following formula (\cite[Lemma $3.2$]{H5}):
\begin{equation}
\label{eq:McMahon} \cF_{A_n}=s_{(n)}\left[\frac{X}{(1-t)(1-q)}\right].
\end{equation}

Denote by $\otimes$ the product on $\Lambda$ induced from $\Phi\circ \chi$ by the
tensor product of representations. We have the following formulas for
$S_n\times T^2$-equivariant $A_n$-modules of finite type:
\begin{proposition}
\label{prop:series} \text{}
\begin{enumerate}
    \item If $0\rightarrow N\rightarrow M\rightarrow P\rightarrow 0$ is an
    exact sequence then $\cF_M=\cF_N+\cF_P$.
    \item $\cF_{M\underset{\IC}{\otimes}N}=\cF_M\otimes \cF_N$.
\end{enumerate}
\end{proposition}

We have the following version of Bott's formula:

\begin{proposition}[Bott's formula] \emph{(\cite[Proposition $3.2$]{H5})} \label{prop:Bott} Let
$F$ be a $T^2$-linearized vector bundle on $\hilb$. Then:
$$
\cF_{\Upsilon(F)}=\sum\limits_{\mu\vdash n} \frac{F_\mu
\cF_{P_n(\xi_\mu)}}{\prod\limits_{x\in
D(\mu)}(1-t^{1+l(x)}q^{-a(x)})(1-t^{-l(x)}q^{1+a(x)})}.
$$
\end{proposition}

The connection between Hilbert schemes and Macdonald polynomials is contained
in the following result:
\begin{theorem}[Haiman]\emph{(\cite[Proposition $3.4$]{H5})} \label{th:Haiman} For $\mu\vdash n$,
$\cF_{P_n(\xi_\mu)}=\tH_\mu$.
\end{theorem}

\section{Action of linearized bundles}

We consider the following problem: each vector bundle $F$ on $\hilb$ induces an
operator on $K(\hilb)$ given by tensor product $G\mapsto F\otimes G$. We look
at a description of the operator one gets on $\Lambda^n$ through
$\Theta:\Lambda^n\rightarrow K(\hilb)$:
$$
\xymatrix{\Lambda^n \ar[rr]^-\Theta\ar@{.>}[d]& & K(\hilb)\ar[d]^{F\otimes - } \\
\Lambda^n \ar[rr]^-\Theta && K(\hilb)}
$$
Since the vector space $K(\hilb)$ has a basis of $T^2$-linearized bundles
(proposition \ref{prop:IN}), it is enough to work with linearized bundles.
Using Haiman's methods (inspired by \cite[Proposition $5.4.9$]{H2}) we get the
following result:

\begin{theorem}
\label{th:action bundle} Let $F$ be a $T^2$-linearized vector bundle on $\hilb$
and denote by $F_\mu\in \IQ[q^{\pm 1},t^{\pm 1}]$ the representation associated to
the action on the fibre at each fixed point. Define an operator
$\nabla_F:\LambdaQ^n\rightarrow \LambdaQ^n$ by
$$
\nabla_F \tH_\mu=F_\mu \tH_\mu \quad \forall \mu\vdash n,
$$
and let $\opE_F:\Lambda^n\rightarrow \Lambda^n$ be
$$
\opE_F:=\left(\omega \nabla_F^*\omega\right)_{q=1,t=1},
$$
where $\nabla^*_F$ is the adjoint operator of $\nabla_F$ for the scalar product
$\langle,\rangle$. Then one has in $K(\hilb)$:
$$
F\otimes \Theta(y)=\Theta(\opE_F(y))\quad \forall y\in \Lambda^n.
$$
\end{theorem}

\begin{proof}
For any $T^2$-linearized vector bundle $G$ on $\hilb$, Bott's formula
\ref{prop:Bott} and Haiman's theorem \ref{th:Haiman} imply:
\begin{align*}
\cF_{\Upsilon(F)}&=\sum_{\mu \vdash n} \frac{F_\mu \tH_\mu}{\prod\limits_{x\in
D(\mu)}(1-t^{1+l(x)}q^{-a(x)})(1-t^{-l(x)}q^{1+a(x)})} ,\\
\cF_{\Upsilon(F\otimes G)}&=\sum_{\mu \vdash n} \frac{F_\mu G_\mu
\tH_\mu}{\prod\limits_{x\in
D(\mu)}(1-t^{1+l(x)}q^{-a(x)})(1-t^{-l(x)}q^{1+a(x)})}.
\end{align*}
By definition of the operator $\nabla_F$ this means that
$\cF_{\Upsilon(F\otimes G)}=\nabla_F \cF_{\Upsilon(G)}$.
The modules (or finite formal sums of modules in $K_{S_n\times T^2}(A_n)$)
$\Upsilon(G)$ and $\Upsilon(F\otimes G)$ admit $S_n\times T^2$-equivariant
finite free resolutions inducing decompositions in the Grothendieck group. This
shows that their Frobenius series satisfy equations:
\begin{align*}
\cF_{\Upsilon(G)}&=\cF_{A_n} \otimes R_G(q,t) ,\\
\cF_{\Upsilon(F\otimes G)}&=\cF_{A_n} \otimes R_{FG}(q,t) ,
\end{align*}
with $R_G(q,t),R_{FG}(q,t) \in \Lambda^n \otimes \IQ[q^{\pm 1},t^{\pm 1}]$. In
the definition of $\Theta$, the Thom isomorphism $\tau$ means that we keep one
fibre and by ``forgetting'' the $T^2$-action we see:
\begin{align*}
\Theta^{-1}(G)&=R_G(1,1) ,\\
\Theta^{-1}(F\otimes G)&= R_{FG}(1,1) .
\end{align*}

So we have to solve the following equation:
$$
\cF_{A_n}\otimes R_{FG}(q,t)= \nabla_F \left(\cF_{A_n} \otimes R_G(q,t)
\right).
$$
With formula (\ref{eq:McMahon}) we get:
$$
s_{(n)}\left[\frac{X}{(1-t)(1-q)} \right]\otimes
R_{FG}(q,t)=\nabla_F\left(\cF_{A_n} \otimes R_G(q,t)\right).
$$
From the formula $p_\lambda \otimes p_\mu=z_\lambda \delta_{\lambda,\mu}
p_\lambda$ (see Manivel \cite{M}) we deduce the plethystic relation:
\begin{align*}
(p_\lambda \otimes p_\mu)\left[\frac{X}{(1-t)(1-q)} \right]&=p_\lambda \otimes
\left(p_\mu
\left[\frac{X}{(1-t)(1-q)} \right] \right) \\
&=\left(p_\lambda \left[\frac{X}{(1-t)(1-q)} \right] \right) \otimes p_\mu.
\end{align*}
Extending by bilinearity and using that $s_{(n)}$ is the unit for the tensor
product, we get the resolution of the equation:
$$
R_{FG}(q,t)=\left( \nabla_F \left(\cF_{A_n} \otimes R_G(q,t)\right) \right)
[X(1-t)(1-q)].
$$
By the duality formula (\ref{eq:duality}), the adjoint operator
$\nabla_F^*$ satisfies \footnote{Observe that the operator $\omega$ commutes
with these plethysms.}:
$$
\nabla_F^*\big(\omega \tH_\mu[X(1-t)(1-q)]\big)=F_\mu \omega
\tH_\mu[X(1-t)(1-q)].
$$
Then,
\begin{align*}
(\nabla_F \tH_\mu)[X(1-t)(1-q)]&=F_\mu \tH_\mu[X(1-t)(1-q)] \\
&= (\omega \nabla_F^* \omega) \left(\tH_\mu[X(1-t)(1-q)] \right).
\end{align*}
Extending this relation by linearity we find:
$$
(\nabla_F f)[X(1-t)(1-q)]=(\omega \nabla_F^* \omega) \left(f[X(1-t)(1-q)]
\right) \quad \forall f.
$$
This applied to $f=\cF_{A_n} \otimes R_G(q,t)$, using again the plethysm
transfer in the tensor product, we get:
\begin{align*}
R_{FG}(q,t)&=(\omega \nabla_F^* \omega) \big( (\cF_{A_n} \otimes R_G(q,t))[X(1-t)(1-q)] \big)  \\
&=(\omega \nabla_F^* \omega)R_G(q,t).
\end{align*}
We evaluate at $t=1,q=1$ (this does work since our operator comes from a
well-defined operator on $\Lambda$) and this gives:
$$
R_{FG}(1,1)=\opE_F R_G(1,1),
$$
hence for any $G$:
$$
\Theta^{-1}(F\otimes G)=\opE_F(\Theta^{-1}(G)),
$$
or equivalently:
$$
F \otimes\Theta(y)=\Theta(\opE_F (y)) \qquad \forall y \in \Lambda^n.
$$
\end{proof}

\begin{remark}
By proposition \ref{prop:IN}, this theorem tells us how the multiplication by a
Schur function $s_\lambda$ acts for the ring structure on $\Lambda$ induced by
$\Theta$: multiplication by $s_\lambda$ means composition with the operator
$\nabla_{P_\lambda^*}$. The weights needed to compute this operator are
$K_{\mu,\lambda}(q^{-1},t^{-1})$.
\end{remark}

\section{The tautological bundle}
\label{s:tautological}

In the case of the tautological bundle $B_n$, we can find a nice formula for
the operator $\opE_{B_n}$. Since the definition of this operator is the same
for all $n$, we shorten the notation and set $\nabla:=\nabla_{B_n}$ and
$\opE:=\opE_{B_n}$ for all $n$, considered as operators defined globally on
$\LambdaQ$ and $\Lambda$ respectively. We state the result as follows:

\begin{theorem}
\label{th:operatorE} Consider the following differential operator on $\Lambda$:
$$
\opE:=\Coeff\left(t^0,\left(\sum\limits_{r\geq 1}
rp_rt^r\right)\exp\left(\sum\limits_{r\geq 1}\frac{\partial}{\partial p_r}
t^{-r}\right)\right).
$$
Then for any $y\in \Lambda$ one has:
$$
B_n\otimes\Theta(y)=\Theta(\opE(y)).
$$
\end{theorem}

\begin{proof}
The proof consists in manipulations of plethystic expressions and some
combinatorics. Observe that the representation $B_n(\xi_\mu)$ is precisely the
polynomial $B_\mu$ defined in \S \ref{ss:Macdonald}, so the operator $\nabla$
is defined by: $\nabla \tH_\mu = B_\mu \tH_\mu$. Since the operator $\Delta$ is characterized by:
$$
\Delta \tH_\mu=\big(1-(1-q)(1-t)B_\mu\big) \tH_\mu,
$$
we have:
$$
\nabla=\frac{1}{(1-q)(1-t)} (id - \Delta),
$$
so by definition of $\Delta$ we get:
$$
\nabla f=\left.\frac{1}{(1-q)(1-t)} \left( f- f\left[X+\frac{(1-t)(1-q)}{z}
\right]\Omega[-zX]\right)\right|_{z^0}.
$$

We first compute a plethystic expression for the operator $\omega \nabla^*
\omega$:

\begin{lemma} The operator $\omega \nabla^* \omega$ is defined by:
$$
(\omega \nabla^* \omega)(f)=\left. \frac{1}{(1-t)(1-q)} \left(f - f \left[
X+\frac{1}{z} \right]\Omega[-zX(1-t)(1-q)]\right) \right|_{z^0}.
$$
\end{lemma}

\begin{proof}[Proof of the lemma]
The duality formula (\ref{eq:duality}) implies:
$$
\nabla^* \big(\omega \tH_\mu[X(1-t)(1-q)]\big)=B_\mu \omega
\tH_\mu[X(1-t)(1-q)].
$$
Setting $S_\mu:=\tH_\mu[X(1-t)(1-q)]$ we get:
\begin{align*}
(\omega \nabla^* \omega) S_\mu&=B_\mu S_\mu \\
&=(B_\mu \tH_\mu)[X(1-t)(1-q)].
\end{align*}
By definition of $\nabla$ we have:
$$
B_\mu\tH_\mu=\left. \frac{1}{(1-t)(1-q)} \left( \tH_\mu -
\tH_\mu\left[X+\frac{(1-t)(1-q)}{z}\right]\Omega[-zX]\right) \right|_{z^0}.
$$
We make the plethystic substitution $[X(1-t)(1-q)]$:
$$
\scriptstyle{B_\mu S_\mu= \left.\frac{1}{(1-t)(1-q)}\left(S_\mu  -
\tH_\mu\left[X+\frac{(1-t)(1-q)}{z}\right][X(1-t)(1-q)]\cdot\Omega[-zX][X(1-t)(1-q)]\right)
\right|_{z^0}} .
$$
By associativity of plethysm we find:
$$
\scriptstyle{B_\mu S_\mu= \left.\frac{1}{(1-t)(1-q)}\left(S_\mu  -
\tH_\mu\left[X(1-t)(1-q)+\frac{(1-t)(1-q)}{z}\right]\cdot\Omega[-zX][X(1-t)(1-q)]\right)
\right|_{z^0}} .
$$
Observe that:
$$
\Omega[-zX][X(1-t)(1-q)]=\Omega[-zX(1-t)(1-q)].
$$
Since $\tH_\mu=S_\mu \left[\frac{X}{(1-t)(1-q)} \right]$ we find:
\begin{align*}
\textstyle{\tH_\mu\left[X(1-t)(1-q)+\frac{(1-t)(1-q)}{z}\right]} &=\textstyle{S_\mu \left[\frac{X}{(1-t)(1-q)} \right]\left[X(1-t)(1-q)+\frac{(1-t)(1-q)}{z}\right] }\\
&=\textstyle{S_\mu \left[ X+\frac{1}{z} \right]}.
\end{align*}
Then,
$$
(\omega \nabla^* \omega) S_\mu=\left. \frac{1}{(1-t)(1-q)} \left(S_\mu - S_\mu
\left[ X+\frac{1}{z} \right]\Omega[-zX(1-t)(1-q)]\right) \right|_{z^0},
$$
hence the result by linearity.
\end{proof}

We deduce a first expression of the operator $\opE$:
\begin{lemma}\label{lemm:opEbase}
For any partition $\lambda=(\lambda_1,\ldots,\lambda_{l(\lambda)})$,
$$
\opE(p_\lambda)=\sum_{\substack{I\subset \{1,\ldots,l(\lambda)\} \\ I \neq
\emptyset}} |\lambda_I| p_{|\lambda_{I}|} p_{\lambda_{\bar{I}}}
$$
where $I$ is a choice of parts, $\lambda_I$ the partition obtained by
preserving these parts, $\lambda_{\bar{I}}$ the complementary partition and
$|\lambda_I|$ the sum of the parts.
\end{lemma}

\begin{proof}[Proof of the lemma]
Starting from the plethystic expression of $\omega\nabla^*\omega$:
$$
(\omega\nabla^*\omega) f=\left. \frac{1}{(1-q)(1-t)} \left(
f-f\left[X+\frac{1}{z}\right]\Omega[-zX(1-q)(1-t)] \right) \right|_{z^0},
$$
we compute in the basis $\{p_\lambda\}$. For a partition
$\lambda=(\lambda_1,\ldots,\lambda_{l(\lambda)})$ we have:
\begin{align*}
\scriptstyle{(\omega\nabla^*\omega)
p_\lambda}&=\scriptstyle{\left.\frac{1}{(1-t)(1-q)} \left( p_\lambda
-\prod\limits_{i=1}^{l(\lambda)} \left(p_{\lambda_i}
+\frac{1}{z^{\lambda_i}}\right) \exp \sum\limits_{r\geq 1}
\frac{-z^r(1-t^r)(1-q^r)}{r}p_r \right) \right|_{z^0}}\\
&=\scriptstyle{ \left.\frac{1}{(1-t)(1-q)} \left( p_\lambda
 - \sum\limits_{I\subset \{1,\ldots,l(\lambda)\}} \frac{\prod\limits_{j\notin
I}p_{\lambda_j}}{\prod\limits_{j \in I}z^{\lambda_j}} \prod\limits_{r \geq 1}
\sum\limits_{m_r \geq 0} \frac{(-1)^{m_r}z^{r
m_r}(1-t^r)^{m_r}(1-q^r)^{m_r}}{m_r!r^{m_r}}p_r^{m_r}\right) \right|_{z^0}}.
\end{align*}
If $I=\emptyset$, we get $p_\lambda$ in the sum, hence this term cancels with
the first $p_\lambda$, so:
$$
(\omega\nabla^*\omega)p_\lambda=- \sum_{\substack{I\subset \{1,\ldots,l(\lambda)\} \\
I \neq \emptyset}} p_{\lambda_{\bar{I}}} \sum_{|\mu|=|\lambda_I|}
(-1)^{l(\mu)}\frac{1}{z_\mu} 
\frac{\prod\limits_{r \geq 1}(1-t^r)^{m_r}(1-q^r)^{m_r}}{(1-t)(1-q)} p_\mu ,
$$
where in the second sum we have set $\mu=(1^{m_1},2^{m_2},\ldots)$. When we
evaluate at $t=q=1$, if the partition $\mu$ has at least two non-zero parts,
the corresponding term cancels. The only remaining case is when $\mu$ is the
partition $(|\lambda_I|)$ and so:
\begin{align*}
\opE p_\lambda&=\sum_{\substack{I\subset \{1,\ldots,l(\lambda)\} \\
I \neq \emptyset}} p_{\lambda_{\bar{I}}} \frac{1}{z_{|\lambda_I|}}
|\lambda_I|^2
p_{|\lambda_I|} \\
&=\sum_{\substack{I\subset \{1,\ldots,l(\lambda)\} \\
I \neq \emptyset}} |\lambda_I|  p_{|\lambda_I|}p_{\lambda_{\bar{I}}}.
\end{align*}
\end{proof}

To finish with, we show that this formula is the same as the formula announced
in the theorem. We start from the end:
$$
\opE=\left. \left(\sum\limits_{r\geq 1}r p_r t^r\right) \exp
\left(\sum\limits_{r\geq 1}\frac{\partial }{\partial p_r}t^{-r}\right) \right|
_{t^0},
$$
and develop the expression:
\begin{align*}
\opE&=\left.\left( \sum_{r \geq 1} r p_r t^r \right) \left( 1+\sum_{k \geq 1}
\frac{1}{k!} \left(
\sum_{r \geq 1} \frac{\partial}{\partial p_r} t^{-r} \right) ^k \right)\right|_{t^0} \\
&=\left.\left( \sum_{r \geq 1} r p_r t^r \right) \left( 1+\sum_{k \geq 1}
\frac{1}{k!} \sum_{n_1,\ldots,n_k \geq 1} \frac{\partial}{\partial p_{n_1}}
\cdots \frac{\partial}{\partial
p_{n_k}} t^{-(n_1+\cdots+n_k)} \right) \right|_{t^0} \\
&=\sum_{k\geq 1} \frac{1}{k!} \sum_{n_1,\ldots,n_k \geq
1}(n_1+\cdots+n_k)p_{n_1+\cdots+n_k} \frac{\partial}{\partial p_{n_1}} \cdots
\frac{\partial}{\partial p_{n_k}}.
\end{align*}

Let
$\lambda=(\lambda_1,\ldots,\lambda_{l(\lambda)})=(1^{\alpha_1},2^{\alpha_2},\ldots)$
be a partition. In the computation of $\opE(p_\lambda)$ with this last formula,
only the $k$-tuples $(n_1,\ldots,n_k)$ consisting of parts of $\lambda$
repeated with a multiplicity less in this $k$-tuple than in $\lambda$
contribute to the sum: for the others, we have $\frac{\partial}{\partial
p_{n_1}} \cdots \frac{\partial}{\partial p_{n_k}}p_\lambda=0$. Then we can
index the remaining $k$-tuples by the choice of $k$ parts in $\lambda$:
$$
I=\{i_1,\ldots,i_k\} \subset \{1,\ldots,l(\lambda)\} \text{ subset of } k
\text{ elements},
$$
and we set $n_1=\lambda_{i_1},\ldots,n_k=\lambda_{i_k}$. We examine in this
way all contributing $k$-tuples, but we have to multiply by $k!$ in order to
permute the parts (the $k$-tuples are not ordered), and divide by the
injectivity defect of the association $I\mapsto \lambda_I$, where $\lambda_I$ is
the sub-partition of $\lambda$ obtained by preserving only the parts selected
by $I$. Setting $\lambda_I=(1^{\alpha^I_1},2^{\alpha^I_2},\ldots)$ one can see
that the injectivity defect is:
$$
\prod_{j\geq 1} \frac{\alpha_j!}{(\alpha_j-\alpha^I_j)!}.
$$
Then,
$$
\opE(p_\lambda)=\sum_{k\geq 1} \sum_{|I|=k}\prod_{j\geq 1}
\frac{(\alpha_j-\alpha^I_j)!}{\alpha_j!}
|\lambda_I|p_{|\lambda_I|}\frac{\partial}{\partial p_{\lambda_{i_1}}} \cdots
\frac{\partial}{\partial p_{\lambda_{i_k}}}p_\lambda.
$$
We have:
$$
\frac{\partial}{\partial p_{\lambda_{i_1}}} \cdots \frac{\partial}{\partial
p_{\lambda_{i_k}}}p_\lambda=\prod_{j \geq 1}
\frac{\partial^{\alpha^I_j}}{\partial p_j^{\alpha^I_j}} p_j^{\alpha_j}=
\prod_{j\geq 1}
\frac{\alpha_j!}{(\alpha_j-\alpha^I_j)!}p_j^{\alpha_j-\alpha^I_j}.
$$
Denote by $\bar{I}$ the parts not selected in $I$. Then:
$$
\frac{\partial}{\partial p_{\lambda_{i_1}}} \cdots \frac{\partial}{\partial
p_{\lambda_{i_k}}}p_\lambda=\prod_{j\geq 1}
\frac{\alpha_j!}{(\alpha_j-\alpha^I_j)!}p_{\lambda_{\bar{I}}},
$$
and so we get the formula we wanted:
$$
\opE(p_\lambda)=\sum_{\substack{I\subset \{1,\ldots,l(\lambda)\} \\ I \neq
\emptyset}} |\lambda_I| p_{|\lambda_{I}|} p_{\lambda_{\bar{I}}}.
$$
\end{proof}

A similar computation for the operator $\opD$ leads to the following formula, to compare with the lemma \ref{lemm:opEbase}:

\begin{lemma}\label{lemm:opDbase}
For any partition $\lambda=(\lambda_1,\ldots,\lambda_{l(\lambda)})$,
$$
\opD(p_\lambda)=\sum_{\substack{I\subset \{1,\ldots,l(\lambda)\} \\ I \neq
\emptyset}} (-1)^{|I|-1}\langle\lambda_I\rangle p_{|\lambda_{I}|} p_{\lambda_{\bar{I}}}
$$
where $\langle \lambda_I\rangle$ is the product of the parts selected by $I$.
\end{lemma}

\section{Interpretations of the similarity between the operators $\opD$ and $\opE$}

\subsection{} We shall first try to get a better understanding of the way the operators $\opD$ and $\opE$ are defined by considering the space $\Lambda$ as the Fock representation of an infinite Heisenberg algebra, in the same spirit as in \cite[\S 2.1]{BZF}. However, the sign conventions in the vertex algebra setup differ from ours.

Let $\kg:=\IQ[t,t^{-1}]$ be the vector space of Laurent polynomials in the indeterminate $t$, considered as a commutative Lie algebra. It has a standard basis of monomials $b_n:=t^n$ for $n\in \IZ$. Let $c:\kg\times\kg\rightarrow \IQ$ be the cocycle
$$
c(f,g):=-\Res_{t=0} f dg.
$$
Equivalently, the cocycle is defined by the rules $c(b_n,b_m)=n\delta_{n+m}$ where $\delta_k$ is the Kronecker symbol, taking the value $1$ if $k=0$ and otherwise $0$.

The \emph{infinite Heisenberg algebra} is the central one-dimensional extension of $\kg$ defined by the cocycle $c$, that is $\kh:=\kg\oplus \IQ \eins$ with basis $\{\eins\}\cup \{b_n\}_{n\in \IZ}$ with the Lie bracket $[\cdot,\cdot]$ characterized by the rules:
$$
[\eins,\cdot]=0,\quad [b_n,b_m]=n\delta_{n+m}\eins.
$$

We consider the vector space $\Lambda=\IQ[p_1,p_2,\ldots]$ as a representation of $\kh$ in the following way:
\begin{itemize}
\item for $k>0$, $b_k$ acts on $\Lambda$ by multiplication by $p_k$ and $b_{-k}$ acts as $-k\frac{\partial}{\partial p_k}$;
\item $b_0$ is the zero morphism;
\item $\eins$ acts as the identity.
\end{itemize}
One checks easily that this defines correctly a representation of the algebra $\kh$ in $\Lambda$. The operators $b_k$ for $k>0$ are called \emph{creation operators} and for $k<0$ \emph{annihilation operators}. Since the creation operators generate
$\Lambda$ from the unit $1\in \Lambda$ and any element in $\Lambda$ can be reduced to the unit $1$ by use of annihilation operators, this is in fact an irreducible representation of $\kh$. Relatively to the \emph{conformal graduation} $\Lambda=\bigoplus\limits_{n=0}^\infty \Lambda^n$, the operators $b_k$ have conformal degree $k$.

Let $\kU:=U(\kh)$ be the universal enveloping algebra of $\kh$. By the
Poincar{\'e}-Birkhoff-Witt theorem, it admits a basis consisting of the unit
$\eins$ and the vectors $b_{\nu_1}\cdots b_{\nu_\ell}$ for all finite
non-increasing sequences $\nu_1\geq \cdots\geq \nu_\ell$ in $\IZ$. The
representation $\kh\rightarrow \End(\Lambda)$ induces an algebra homomorphism
$\kU\rightarrow \End(\Lambda)$ allowing to interpret all polynomial
differential operators on $\Lambda$ as elements of $\kU$. In order to extend
this interpretation to the operators $\opD$ and $\opE$, we consider the
following completion of the algebra $\kU$. Define the left ideals generated by
annihilation operators:
$$
B_k:=\kU\cdot \left\{ b_{-\nu_1}\cdots b_{-\nu_\ell} ,|\, \nu_1\geq \cdots\geq \nu_\ell>0 \text{ and } \sum\limits_{i=1}^\ell \nu_i\geq k \right\},\quad k\in \IZ.
$$
Taking the $B_k$'s as a fundamental system of neighborhoods at the origin, one gets a topology on $\kU$ whose corresponding completion is the algebra defined as the inverse limit of the inverse system $\theta_{k+1}:\kU/B_{k+1}\rightarrow \kU/B_{k}$:
$$
\hat{\kU}:=\varprojlim \kU/B_{k}.
$$

\begin{lemma}
The action $\kU\rightarrow \End(\Lambda)$ extends to an action $\hat{\kU}\rightarrow \End(\Lambda)$.
\end{lemma}

\begin{proof}
Let $\xi\in\hat{\kU}$ given by a coherent system $(\xi_k)_{k\in \IZ}$ where
$\xi_k\in \kU/B_k$ with the condition $\theta_{k+1}(\xi_{k+1})=\xi_k$ for all
$k$. In order to define the action of $\xi$ on an element $v\in \Lambda$, one
can assume that $v\in \Lambda^n$ for some $n$. Then one remarks that for all
$k>n$, $B_k\cdot v=0$ since all elements of $B_k$ begin with an annihilation of
order at least $k$. Hence, if $\phi_k\in \kU$ represents the class $\xi_k$, the
evaluation $\phi_k(v)$ only depends on $\xi_k$ when $k>n$. The condition
$\theta_{k+1}(\xi_{k+1})=\xi_k$ also implies that this evaluation does not
depend on the choice of the integer $k>n$. This defines correctly the action
$\xi(v)$ and it is straightforward to check that it is an algebra action.
\end{proof}

In the proof of the theorem \ref{th:operatorE}, we wrote the differential operator $\opE$ as:
$$
\opE=\sum_{k\geq 1} \frac{1}{k!} \sum_{n_1,\ldots,n_k \geq
1}(n_1+\cdots+n_k)p_{n_1+\cdots+n_k} \frac{\partial}{\partial p_{n_1}} \cdots
\frac{\partial}{\partial p_{n_k}}.
$$
We can then interpret the operator $\opE$ as a well-defined element of the completion $\hat{\kU}$ given by:
$$
\opE=\sum_{\ell\geq 1} \frac{(-1)^\ell}{\ell!} \sum_{\nu_1,\ldots,\nu_\ell \geq
1}\frac{\nu_1+\cdots+\nu_\ell}{\nu_1\cdots\nu_\ell}b_{\nu_1+\cdots+\nu_\ell} b_{-\nu_1} \cdots
b_{-\nu_\ell}.
$$
A similar computation gives for the operator $\opD$ the following differential expression:
$$
\opD=\sum_{k\geq 1} \frac{(-1)^{k-1}}{k!} \sum_{n_1,\ldots,n_k \geq
1}(n_1\cdots n_k)p_{n_1+\cdots+n_k} \frac{\partial}{\partial p_{n_1}} \cdots
\frac{\partial}{\partial p_{n_k}},
$$
so the corresponding element in $\hat{\kU}$ is:
$$
\opD=-\sum_{\ell\geq 1} \frac{1}{\ell!} \sum_{\nu_1,\ldots,\nu_\ell \geq
1}b_{\nu_1+\cdots+\nu_\ell} b_{-\nu_1} \cdots b_{-\nu_\ell}.
$$

Now that our operators are well-defined, we can investigate their similarity.

\subsection{First interpretation}
Let $\Omega:\kh\rightarrow \kh$ be the linear isomorphism defined by $\Omega(\eins)=\eins$ and for $k>0$ by $\Omega(b_k)=-kb_k$ and $\Omega(b_{-k})=-\frac{1}{k}b_{-k}$. In fact, this is a Lie algebra isomorphism since $\Omega$ is compatible with the Lie bracket. This extends naturally to an algebra isomorphism of the universal enveloping algebra, that we also denote by $\Omega:\kU\xrightarrow{\sim} \kU$. Since $\Omega^{-1}(B_k)=B_k$ for all $k$, this map is continuous for the chosen topology and so induces an isomorphism of completed algebras $\hat{\Omega}:\hat{\kU}\xrightarrow{\sim} \hat{\kU}$. The similarity between the operators $\opD$ and $\opE$ is then nothing else than the following observation:
$$
\hat{\Omega}(\opD)=\opE.
$$

\subsection{Second interpretation}
Let $\Pi:\Lambda\rightarrow \Lambda$ be the algebra isomorphism characterized
by $\Pi(p_r)=-\frac{1}{r}p_r$, so that for a partition
$\lambda=(\lambda_1,\ldots,\lambda_k)$ one has
$\Pi(p_\lambda)=(-1)^{l(\lambda)}\frac{1}{\langle \lambda\rangle}p_\lambda$
where $\langle \lambda\rangle:=\prod\limits_{i=1}^{k}\lambda_i$. Then one has
the following conjugation property:

\begin{proposition} \label{prop:conjbyPi}The operators $\opD$ and $\opE$ are conjugated by the isomorphism $\Pi$:
$$
\opD \circ \Pi=\Pi \circ \opE.
$$
\end{proposition}

\begin{proof}
From the expression of the operators $\opD$  and $\opE$ given in the lemmas
\ref{lemm:opDbase} and \ref{lemm:opEbase} if follows, since $\Pi$ is an algebra
homomorphism, that for any partition $\lambda$ one has:
\begin{align*}
\Pi(\opE(p_\lambda))&=\sum_{\substack{I\subset \{1,\ldots,l(\lambda)\} \\ I
\neq \emptyset}} |\lambda_I|
\left(-\frac{1}{|\lambda_I|}p_{|\lambda_{I}|}\right)
\left((-1)^{l(\lambda)-|I|}\frac{1}{\langle \lambda_{\bar{I}}\rangle} p_{\lambda_{\bar{I}}}\right) \\
&= \opD(\Pi(p_\lambda)).
\end{align*}
\end{proof}

\subsection{Relation with the Chern Character}

We consider the following situation (the restrictions of the operator $\opD$ and $\opE$ to $\Lambda^n$ are also denoted by $\opD$ and $\opE$):
$$
\xymatrix{ \Lambda^n \ar@(l,u)[]|{\opD} \ar[rr]^-\Psi && H^*(\hilb) \\
\Lambda^n \ar@(d,l)[]|{\opE} \ar@{.>}[u]^\Gamma \ar[rr]^-\Theta &&  K(\hilb)
\ar[u]^{ch} }
$$
The bundle $B_n$ acts on $K(\hilb)$ by tensor product whereas $ch(B_n)$ acts on
$H^*(\hilb)$ by cup product. The operators $\opD$ et $\opE$ play similar roles
in the understanding of the product induced by the cup product via $\Psi$ and
the study of the product induced by the tensor product via $\Theta$.

We remark the following property:
\begin{proposition} \label{prop:conjbyGamma}The operators $\opD$ and $\opE$ are conjugated by the morphism $\Gamma$:
$$
\Gamma \circ \opE=\opD \circ \Gamma.
$$
\end{proposition}

\begin{proof}
By theorem \ref{th:operatorE} we have:
$$
\Theta(\opE (y))=B_n \otimes \Theta (y) \qquad \forall y\in \Lambda^n,
$$
and since the Chern character is an algebra morphism:
$$
(ch\circ\Theta)(\opE (y))=ch(B_n) \cup (ch\circ\Theta) (y).
$$
Since $ch\circ \Theta=\Psi \circ \Gamma$ we get with theorem
\ref{th:operatorD}:
\begin{align*}
(\Psi \circ \Gamma)(\opE (y))&=ch(B_n) \cup (\Psi \circ \Gamma) (y) \\
&=\Psi(\opD (\Gamma y)).
\end{align*}
Since $\Psi$ is an isomorphism, this means $\Gamma \circ \opE=\opD \circ
\Gamma$.
\end{proof}

The composite morphism $\Gamma$ has been studied in details in \cite{SB2,SB1},
where we found the following explicit formula:

\begin{theorem}[Boissi{\`e}re]{\cite[Theorem $5.2$]{SB2}}
\label{th:operatorGamma} For each partition $\mu$ of $n$,
$$
\Gamma(p_\mu)=\sum_{\nu\vdash n} z_\nu^{-1} \sum_{\lambda\vdash n}
\frac{1}{h(\lambda)} \chi^{\lambda}_\nu \chi^{\lambda}_\mu
\Coeff\left(t^{n-l(\nu)},P_{\lambda,\mu}(t)\right) p_\nu,
$$
where the power series $P_{\lambda,\mu}(t)$ is given by Taylor expansion of the
expression:
$$
P_{\lambda,\mu}(t)=\frac{e^{-n(\lambda)t}\prod\limits_{x\in
D(\lambda)}(1-e^{h(x)t})}{\prod\limits_{i=1}^{l(\mu)}(1-e^{\mu_it})}.
$$
In particular, the operator $\Gamma$ has the following shape:
$$
\Gamma(p_\mu)=\frac{(-1)^{n-l(\mu)}}{\langle
\mu\rangle}p_\mu+\sum_{\substack{\nu\vdash n\\l(\nu)<l(\mu)}} g_{\mu,\nu}p_\nu
$$
for some coefficients $g_{\mu,\nu}$.
\end{theorem}

\section{Filtration by the age}

Let us fix an integer $n$ and for any partition $\lambda$ of $n$, set
$\age(\lambda):=n-l(\lambda)$. This is exactly half the cohomological degree in
$H^*(\hilb)$ through the map $\Psi$, and the terminology ''age" comes from the
relation with the orbifold degree shifting number (see \cite{SB1} and
references therein). Define the following subspaces of $\Lambda^n$:
$$
V_i:=\IQ\left\{ p_\lambda \,|\, \age(\lambda)= i\right\}, \quad
F_i:=\IQ\left\{ p_\lambda \,|\, \age(\lambda)\geq i\right\},\quad i=0,\ldots,n-1.
$$
This constructs a decreasing flag
$$
\Lambda^n=F_0\supset F_{1}\supset\cdots\supset F_{n-1} \supset F_n=\{0\}
$$
with factors $F_{i}/F_{i+1}\cong V_i$, and also a graded decomposition
$\Lambda^n=\bigoplus_{i=0}^{n-1} V_i$.

Denote by $P^{\age}$ the parabolic subgroup of the linear group
$\GL(\Lambda^n)$ consisting of automorphisms stabilizing this flag, that is all
$f\in \GL(\Lambda^n)$ so that $f(F_i)\subset F_{i}$ for all $i$. In the graded
decomposition of $\Lambda^n$, such an $f$ is given by a lower-triangular
matrix. For such an automorphism $f\in P^{\age}$, denote by $f_k$, $k=0,..,n-1$
the graded component of $f$ so that $f_k(V_i)\subset V_{i+k}$ for all $i$. In
the graded decomposition, $f_k$ is the submatrix of $f$ given by the
$(i+k,i)$-blocks. One can then decompose $f=f_0+f_1+\cdots+f_{n-1}$.

\begin{lemma}\text{}
\begin{enumerate}
\item The operators $\opD,\opE$ and $\Gamma$ lie in $P^{\age}$.

\item $\opE_0=\opD_0=n\cdot\Id$ and $\Gamma_0=(-1)^n\Pi$.
\end{enumerate}
\end{lemma}

\begin{proof}
This is obvious with the formulas given in lemma \ref{lemm:opEbase},
lemma \ref{lemm:opDbase} and the theorem \ref{th:operatorGamma}.
\end{proof}

\section{Observations on the induced ring structure}

In this section and the next one, we conclude this study by some observations
on the ring structure on $\Lambda^n$ induced from the tensor product of sheaves
on $K(\hilb)$ by the Bridgeland-King-Reid map $\Theta:\Lambda^n \rightarrow
K(\hilb)$. We denote by $\odot$ the product so created on $\Lambda^n$.

\subsection{Basic facts}

Denote the decreasing \emph{topological} filtration of $K(\hilb)$ defined by
the codimension of the support of coherent sheaves by
$$
F_dK(\hilb):=\IQ\left\{\cF\,|\,\codim \Supp \cF\geq d\right\}.
$$
The ring structure on $K(\hilb)$ is compatible with this filtration (see
\cite{BS,FL}), that is $F_dK(\hilb)\cdot F_{d'}K(\hilb)\subset
F_{d+d'}K(\hilb)$. The Chern character $ch:K(\hilb)\rightarrow H^*(\hilb)$ is a
ring isomorphism and is compatible with this filtration. As explained in
\cite[Theorem $5.4$]{SB2}, the shape of the map $\Gamma$ shows that the
morphism $\Theta$ is compatible with the topological filtration of $K(\hilb)$
and the age filtration of $\Lambda$: with our notation, $(\Gamma\circ\Theta^{-1})(F_dK(\hilb))=F_d$. The diagram
$$
\xymatrix{ \Lambda^n \ar[rr]^-\Psi && H^*(\hilb) \\
\Lambda^n \ar[u]^\Gamma \ar[rr]^-\Theta &&  K(\hilb) \ar[u]^{ch} }
$$
thus induces a graded diagram with respect to the various filtrations:
$$
\xymatrix{ \gr\Lambda^n \ar[rr]^-\Psi && H^*(\hilb) \\
\gr\Lambda^n \ar[u]^{\Gamma_0} \ar[rr]^-{\gr\Theta} &&  \gr K(\hilb)
\ar[u]^{\gr ch} }
$$
where $\Psi$ and $\gr ch$ are isomorphisms of graded rings. All this means the
following for the ring structure $(\Lambda^n,\odot)$:
\begin{itemize}
\item the unit is the Schur function $s_{(n)}$ corresponding to the trivial
representation of $S_n$;

 \item the product is filtered with respect to the age
filtration: $F_d\odot F_{d'}\subset F_{d+d'}$;

 \item the induced graded product on
$(\gr\Lambda^n,\gr\odot)$ is isomorphic to the graded ring
$(\gr\Lambda^n,\cup)$ induced from the cup product on the cohomology by the
isomorphism $\Psi$ and explicitely described in \cite{LS}: the isomorphism is
just the rescaling given by $\Gamma_0=(-1)^n\Pi$.
\end{itemize}

\subsection{Action of the Adams powers}

For a quasi-projective algebraic variety $X$ and any integer $j\geq 0$, the
Adams operator $\psi^j:K(X)\rightarrow K(X)$ is the ring homomorphism
characterized by the assumption that $\psi^j$ is functorial with respect to $X$
and for a line bundle $L$ over $X$ one has $\psi^j(L)=L^j$. This determines
uniquely $\psi^j$ by splitting construction (see \cite{FL}). Let $E$ be a
vector bundle on $X$ an denote the homogeneous components of its Chern
character by $ch E=\sum_{k\geq 0}ch_k(E)$ with $ch_k(E)\in H^{2k}(X)$. Then one
checks that the Chern character of $\psi^j(E)$ is:
$$
ch(\psi^j(E))=\sum_{k\geq 0}j^k ch_k(E).
$$

In the case of the Hilbert scheme of points on the affine plane $\hilb$, since
the components of the Chern character of the tautological bundle $B_n$ generate
the cohomology ring $H^*(\hilb)$, this means that the Adams powers $\psi^j B_n$
for $j=0,1,\ldots,n-1$ form a set of generators of the ring $K(\hilb)$, since
the Vandermonde matrix $\Big(j^k\Big)_{j,k=0,\ldots,n-1}$ is invertible.

This implies that the operators $\opE_{\psi^jB_n}$ describing the action of
$\psi^jB_n$ by tensor product, considered as operators on $\Lambda^n$ by the McKay correspondence
$\Theta$, give a complete description of the product $\odot$. One has
$\opE_{{\psi^0B_n}}=\Id$, and $\opE_{\psi^1B_n}=\opE$ has been completely
described. In order to compute the other operators with theorem
\ref{th:action bundle}, one has to use the easy fact that the polynomials giving the
action at a fixed point are $B_\mu(q^j,t^j)=p_j[B_\mu]$. Numerical evidence
suggests the following conjecture:

\begin{conjecture} The operator $\opE_{\psi^jB_n}$ describing the action of an
Adams power $\psi^jB_n$ in $\Lambda^n$ is given by the formula:
$$
\opE_{\psi^jB_n}=\sum_{k\geq 0}j^k \opE_k.
$$
\end{conjecture}

Note that it is clear that, at the cohomological level, the action of
$ch(\psi^j(B_n))$ by cup product, viewed on $\Lambda^n$ by $\Psi$, is given by
$\sum_{k\geq 0}j^k \opD_k$, since the map $\Psi$ is graded. This does not apply
to the K-theoretical setup since the map $\Theta$ is only filtered.

\section{Numerical tables}

One can in fact compute completely the table of the product $\odot$ in the
basis of Newton functions and Schur functions by making use of the explicit
formula for the operator $\Gamma$ given in theorem \ref{th:operatorGamma} and an
implementation of the cup product $(\gr\Lambda^n,\cup)$: the ring
$(\Lambda^n,\odot)$ is indeed given by transfer from the cup product. This can
be achieved with the software SINGULAR \cite{SING}. As an illustration of our
observations, we give here the first tables:
\begin{itemize}
\item $n=2$:

$$
\begin{array}{|c||c|c|}
  \hline
    \odot & p_{1,1} & p_{2} \\
  \hline\hline
  p_{1,1} & 2p_{1,1}-2p_2 & 2p_2 \\
  \hline
  p_{2} & 2p_2 & 0 \\
  \hline
\end{array}
$$

$$
\begin{array}{|c||c|c|}
  \hline
    \odot & s_{2} & s_{1,1} \\
  \hline\hline
  s_{2} & s_{2} &  s_{1,1} \\
  \hline
  s_{1,1} & s_{1,1} & -s_2+2s_{1,1} \\
  \hline
\end{array}
$$

\item $n=3$:
$$
\begin{array}{|c||c|c|c|}
\hline
  \odot & p_{1,1,1} & p_{2,1} & p_{3} \\
  \hline\hline
  p_{1,1,1} & 6p_{1,1,1}-18p_{2,1}+15p_3 & 6p_{2,1}-9p_3 & 6p_3 \\
  \hline
   p_{2,1} & 6p_{2,1}-9p_3 & 3p_3 & 0 \\
   \hline
  p_{3} & 6p_3 & 0 & 0 \\
  \hline
\end{array}
$$

$$
\begin{array}{|c||c|c|c|}
\hline
  \odot & s_{3} & s_{2,1} & s_{1,1,1} \\
  \hline\hline
  s_{3} & s_{3} & s_{2,1} & s_{1,1,1} \\
  \hline
   s_{2,1} & s_{2,1} & -s_3+s_{2,1}+3s_{1,1,1} & s_{3}-2s_{2,1}+5s_{1,1,1} \\
   \hline
  s_{1,1,1} & s_{1,1,1} & s_{3}-2s_{2,1}+5s_{1,1,1} & 2s_3-3s_{2,1}+5s_{1,1,1} \\
  \hline
\end{array}
$$

\newpage

\item $n=4$:
\begin{center}
\begin{rotate}{-90}
$
\begin{array}{|c||c|c|c|c|c|}
\hline
  \odot & p_{1,1,1,1} & p_{2,1,1} & p_{2,2} & p_{3,1} & p_4 \\
  \hline\hline
  p_{1,1,1,1} & 24p_{1,1,1,1}-144p_{2,1,1}+72p_{2,2}+240p_{3,1}-240p_4 & 24p_{2,1,1}-24p_{2,2}-72p_{3,1}+104p_4 & 24p_{2,2}-48p_4 & 24p_{3,1}-48p_4 & 24p_4 \\
  \hline
  p_{2,1,1} & 24p_{2,1,1}-24p_{2,2}-72p_{3,1}+104p_4 & 4p_{2,2}+12p_{3,1}-32p_4 & 8p_4 & 8p_4 & 0 \\
  \hline
  p_{2,2} &  24p_{2,2}-48p_4& 8p_4 & 8p_4 & 0 & 0 \\
  \hline
  p_{3,1} &  24p_{3,1}-48p_4 & 8p_4 & 0 & 0 & 0 \\
  \hline
  p_4 & 24p_4 & 0 & 0 & 0 & 0 \\
  \hline
\end{array}
$
\end{rotate}
\hspace{6cm}
\begin{rotate}{-90}
$
\begin{array}{|c||c|c|c|c|c|}
\hline
  \odot & s_{4} & s_{3,1} & s_{2,2} & s_{2,1,1} & s_{1,1,1,1} \\
  \hline\hline
  s_{4} & s_4 & s_{3,1} & s_{2,2} & s_{2,1,1} & s_{1,1,1,1} \\
  \hline
  s_{3,1} & s_{3,1} & \begin{array}{c}-s_{4}+s_{3,1}\\-s_{2,2}+3s_{2,1,1}\end{array} & \begin{array}{c}-s_{3,1}+2s_{2,2}\\+s_{2,1,1}+2s_{1,1,1,1}\end{array} & \begin{array}{c}s_4-s_{3,1}-2s_{2,2}\\+3s_{2,1,1}+6s_{1,1,1,1}\end{array} & \begin{array}{c}-s_4+2s_{3,1}-s_{2,2}\\-3s_{2,1,1}+9s_{1,1,1,1}\end{array} \\
  \hline
  s_{2,2} & s_{2,2} & \begin{array}{c}-s_{3,1}+2s_{2,2}\\+s_{2,1,1}+2s_{1,1,1,1}\end{array} & \begin{array}{c}-s_{3,1}+3s_{2,2}\\-s_{2,1,1}+4s_{1,1,1,1}\end{array} & \begin{array}{c}s_{2,2}-2s_{2,1,1}\\+10s_{1,1,1,1}\end{array} & \begin{array}{c}-2s_4+3s_{3,1}\\-5s_{2,1,1}+10s_{1,1,1,1}\end{array} \\
  \hline
  s_{2,1,1} & s_{2,1,1} & \begin{array}{c}s_4-s_{3,1}-2s_{2,2}\\+3s_{2,1,1}+6s_{1,1,1,1}\end{array} & \begin{array}{c}s_{2,2}-2s_{2,1,1}\\+10s_{1,1,1,1}\end{array} & \begin{array}{c}-s_4+4s_{3,1}-4s_{2,2}\\-6s_{2,1,1}+24s_{1,1,1,1}\end{array} & \begin{array}{c}-5s_4+8s_{3,1}-2s_{2,2}\\-11s_{2,1,1}+21s_{1,1,1,1}\end{array} \\
  \hline
  s_{1,1,1,1} & s_{1,1,1,1} & \begin{array}{c}-s_4+2s_{3,1}-s_{2,2}\\-3s_{2,1,1}+9s_{1,1,1,1}\end{array} & \begin{array}{c}-2s_4+3s_{3,1}\\-5s_{2,1,1}+10s_{1,1,1,1}\end{array} & \begin{array}{c}-5s_4+8s_{3,1}-2s_{2,2}\\-11s_{2,1,1}+21s_{1,1,1,1}\end{array} & \begin{array}{c}-5s_4+7s_{3,1}-s_{2,2}\\-9s_{2,1,1}+14s_{1,1,1,1}\end{array} \\
  \hline
\end{array}
$
\end{rotate}
\end{center}
\end{itemize}

\newpage

\nocite{*}
\bibliographystyle{amsplain}
\bibliography{BiblioMultBKRHilbert}

\end{document}